\documentclass{article}
\usepackage{amsmath}
\usepackage{amssymb}
\newcommand{\di}{\displaystyle}
\newcommand{\Z}{\mathbb{Z}}
\title{Group ring cryptography}
\date{}
\author{Barry Hurley \& Ted Hurley\\ National University of Ireland Galway.}
\begin{document}
\maketitle


\begin{abstract} 
 Cryptographic systems are derived using units in  group rings. 
Combinations of 
types of units in group rings give  units not of any particular type.  
This includes cases of taking  powers of  
 units and products of such powers and  
adds the complexity of the {\em discrete logarithm} problem to the system.

The method enables encryption and  (error-correcting) coding  
to be combined within one system. 


These group ring cryptographic systems may  
be combined in a neat way 
 with existing
 cryptographic systems, such as
 RSA, and a combination has the combined strength of both systems.

Examples are given. 
 
\end{abstract}


\section{Introduction}
Cryptography is the key element in electronic security systems and has
many uses. Included in cryptography are public key encryption, data
encryption standards, key exchange systems, digital signatures and others.
Probably the best known cryptographic systems  are the RSA system, the Diffie
Hellmann key exchange system and the discrete logarithm methods. 
In recent years there has been much 
interest in the use of {\em Elliptic Curves}, see \cite{koblitz,koblitz1}.

\footnotetext{{\bf Keywords:} Group ring, units,
  cryptography,  public key
  cryptography, coding.} 

These known systems work
within commutative groups.
The RSA system works  in
$U(\mathbb{Z}_n)$, the units of $\Z_n$, the integers modulo $n$,  
where $n$ is the product of two large primes. In this system the  order of
$U(\mathbb{Z}_n)$ is kept secret. Other systems are based on the units
or multiplicative group 
of a finite field $GF(q), \, q = p^n$ with $p$ a prime.

The security of RSA depends on the difficulty of factoring primes; 
others depend on the difficulty of the
{\it discrete logarithm problem}.


There has also been some interest in constructing  public key
cryptosystems via combinatorial group theory. This can be considered 
as a possible generalisation of the discrete logarithm problem to arbitrary
groups using  the difficulty of the so-called {\em conjugacy search
problem}. See for example \cite{ans},  \cite{cryp},
\cite{attack}, \cite{deh} and the references therein; note that
\cite{attack} discusses attacks on those that use Braid
groups and looks at other possibilities. 


The technique here is to use appropriate elements of the group
of units,
$U(RG)$, of a group ring $RG$ for encryption and decryption. 
Many construction
 methods work in many different group rings giving rise to 
 different unique systems of public key encryption.



As group rings may also be used for designing codes, 
the method enables  encryption and (error-correcting) 
coding to be integrated within the same group ring
 system.
 This has potential in terms of complexity reduction, cost savings in terms of
chip design, and the number of possible applications that
could 
benefit from cheap, secure and reliable communication.



In the group ring method, units of particular types may be combined to
  give a unit which is not of any of the constituent types nor of any
  known type. It is 
  computationally 
impossible to obtain the original units from the product without
  information on the units involved. Compare this to RSA in which the two
  primes  are combined to give an element/number from which it is
  impossible computationally to retrieve the original primes. 
To decrypt
      in the group ring method, information on the individual units (primes!)
  involved must be known  {\em as well as} information on the inverses of
  each of these components. A component may itself be a power
  of an unknown unit adding the difficulty of the {\em
  discrete logarithm problem} to the system. 

The group ring method may if required be combined  
with existing systems, such as  RSA, in a
particularly neat way giving a system which has the combined strength of
both systems. 

When the keys have been constructed the time for
encryption and decryption is of order  $n\log n$ at most.

Illustrative examples are given which demonstrate the techniques. 
These examples are not exhaustive and  the lengths
and numbers are kept relatively small here so they can be displayed.

\section{Group ring cryptography}

The {\it group ring} $RG$ of a group $G$ over a  ring $R$ 
is the set of all formal sums 

\[\sum_{i=1}^{n} \alpha_ig_i\]  with $\alpha_i \in R$ and $g_i \in
G$ 

where only a finite
 number of the coefficients (elements of $R$)
are non-zero. See \cite{sehgal} for further information on group rings. 
We assume that the ring $R$ has an identity element but there is no
further restriction on the ring used; in particular it does not
necessarily have to be a field. A group ring is a ring
itself and it's possible to choose the ring $R$ in the group ring $RG$
to be a group ring itself and this can sometimes be useful. 

The set of invertible elements in a ring $H$ is a group, $U(H)$, called
the {\it the group of units of $H$}. 
A zero-divisor in a ring $H$ is an element $u\not = 0
\in H$ such that there is an element $v\in H, v\not = 0$ with $uv=0$.
There is extensive literature on units and zero-divisors in group
rings.

If $u \in U(RG)$ then the inverse of $u$ is denoted by
$u^{-1}$ and satisfies $u*u^{-1} = e = u^{-1}*u$ where $e$ is the
identity element of $RG$ and $*$ denotes multiplication in $RG$. 

\subsection{Method}
 The method of encryption is performed  as follows. 

\begin{itemize}
\item Suppose first of all the data
 or number 
 to be encrypted, $w$ say,  is expressible with digits in a  
 base and whose number of digits in that base 
is less than or equal to  the order of $G$. 

 
The digits of $w$ are considered as
elements of the ring $R$. It is always possible to find such
  a ring as for example the integers $\Z$ will suffice. In many cases
  $R$ will  be the integers $\Z$ or the integers modulo $q$, $\Z_q$, 
  but elements over finite fields or other systems can also be used.

Then $w = \beta_1 \beta_2\ldots\beta_{n} $ with
 the $\beta_i \in R$. Hence $w$ can be uniquely represented as a group
      ring element $\hat{w} = 
\di\sum_{i=1}^{n}\beta_ig_i$ with the $\beta_i$ as the digits, which
 are  elements of $R$, and the $G = \{g_1, g_2, \ldots, g_n\}$ is a
 listing of the elements of $G$.

Thus given a unit $u = \di\sum_{i=1}^{n}\alpha_ig_i$  and an element $w
= \displaystyle{\sum_{i=1}^{n}\beta_ig_i}$  the encryption is given by $w
\mapsto w * u$, where $*$ signifies multiplication in the group
ring. An
encryption may also be given by $w
\mapsto u * w$ and this is different in general.   

We can also add extra digits, and these can be used in various ways
such as for error correcting codes.

The encryption depends on the listing of  the elements of $G$ and a
further complexity could be introduced by permuting the order.    

\item If the number of digits in $w$ is greater than or equal to the order
 of $G$
 then it may
be encrypted by breaking it into smaller blocks in each of which
the number of digits is less than or equal to
 the order of $G$.
Each block is then encrypted.

  \item In the same way an encryption (block encryption) may then also
be applied to the blocks; decrypting can only be performed by
knowing both keys, the key for decrypting  each block and the key  for
decrypting the  block encryption.

One way of applying block encryption is as follows: Let the blocks be
$B_1,B_2,\ldots,B_t$. Consider then $ w =  \di\sum_{i=1}^qB_ih_i$ with the
$h_i$ as the elements of another, or the same,  group $H$ whose order
is greater than
$t$. Assume $B_i = 0$ for $i
> t$ and the $B_i \in R$ for some ring $R$. ($R$ could for example be $\Z$.) 
Then encode by $ w \mapsto w*u$ where $u$ is a unit in $RH$.
\end{itemize}

In some applications it is true that $w*u = u*w$ but this is not
always the case. If $w*u = u*w$ for all $w$ then the
encryption/decryption is said to be {\em commutative}; otherwise it is
said  to be {\em
non-commutative}. Indeed non-commutativity can add another dimension
of complexity.
When non-commutative the encryption $w \mapsto u*w$ is different to 
$w \mapsto w*u$. Further an encrytion could be performed by 
$w\mapsto v*w*u$ for units $v,u$ which in the non-commutative case 
is not the same as $w\mapsto v*u*w$.

The inverse, $u^{-1}$, of $u$ must be known in order to perform the
decryption. The decryption is  $c \mapsto
  u^{-1}*c$ when the encryption is $w \mapsto u*w$; the decryption
  is   is $c\mapsto c*u^{-1}$ when the encryption is $w \mapsto w*u$
  and similarly for others. The inverse may be the  product
  of inverses of certain units. 


\subsection{Large numbers.} When working out a unit to be used for
encryption or decryption it is sometimes the case that large numbers
will occur. We can get over this problem as follows.
Let $u$ be a unit which involves large numbers and which is to be
used for encryption. Suppose the data to be encrypted 
only involves non-negative integers less than some number $m$. Let this
data be $q = \di\sum_{i=1}^{n}\gamma_ig_i$ with $ 0 \leq \gamma_i <
m$. Suppose the encryption is given by $q \mapsto q * u $ and the
decryption is $f \mapsto f * v$ where $v$ is the inverse of $u$.
Then in integer coefficients:

$$q = q * u * v$$

Let $u_{m}$ denote $u$ with its coefficients taken  modulo $m$
and similarly  $v_{m}$ is $v$ with coefficients modulo $m$. Then
we see that

$$q = q * u_{m} * v_{m}$$

from the condition on the coefficients of $q$.

Thus the encryption can be given by  $q \mapsto q * u_{m} $
and the decryption is then $f \mapsto f * v_{m}$ where  the
 the coefficients here are worked out modulo $m$.

Note that $m$ can be taken as large as the computer system and
software will allow. We can thus in this situation work in $\Z_m$ the
integers modulo $m$ or the finite field $GF(m)$ when $m$ is a prime or
a power of a prime. 

\subsection{Products and powers of units.} The product of units is also
a unit. Different {\em types} of units are  known; further information may be
obtained in \cite{sehgal}
and in the references therein. 
The product of a  particular type of units is in most case not of this
type or indeed of any known type. Thus the encryption can be further
made more difficult to decrypt by taking the product of two or more
units and using these to encrypt and decrypt. The product in expanded
form only is given. Compare this to RSA in which two large primes are
multiplied together but only the product, which is not a prime, 
 is given or known to the public. 

For example 
given two units $q, t$  and an element $w$ to be encrypted, the
encryption is given by  $w
\mapsto w*(q*t)$ and the decryption is given by $p \mapsto
p*t^{-1}*q^{-1}$ --  note the order of inversion is important for
non-commutative units. Now $q*t$ only is given and not $q$, $t$
individually but the decoding is done step-by-step. The product will
almost always not be of a known type even when the types of $q,t$ are
known. Similarly we can use a product of a number of units. 

Taking powers of a particular type of unit for which
the type is  known will not be of this type and 
introduces the {\em difficulty of the discrete logarithm problem},
i.e. knowing the power of the unit it is not computationally possible
to deduce the unit.

\subsection{Combining with existing systems.} The group ring unit method of
encryption may also be combined with known methods such as RSA and this
then gives an encryption method which is stronger than the RSA and unit
group ring method. Examples of this are given below.

In RSA the integer $n$ is chosen as the product of two (large) primes, the
alphabet is $\Z_N$ for some positive integer $N$ and $k = [\log_Nn]$,
the integer part of $\log_Nn$. A word $m$ to be encrypted has
the form $m_1\ldots
m_k$ corresponding to the integer $\di\sum_{i=1}^{k}m_iN^{k-i}$. The
encryption is performed by
$m \mapsto m^e \mod n$ and the ciphertext $c = m^e \mod n$ is written
as $ c = \di\sum_{i=0}^{k}c_iN^{k-i}$. 

At any stage either before or after this encryption, or both before and
after, a unit encryption
may be introduced. Consider then $m = \di\sum_{i=1}^{k}m_ig_i$ and a unit
$u = \di\sum_{i=1}^{t}\alpha_ig_i$ where $t \geq k$. First of all encrypt $m$
by $m \mapsto m*u = m^{'}, \, \text{say}$. The unit used is the given
unit of the intended receiver. Then $m^{'}$ is encrypted in the normal
RSA way to give a $c^{'}$. At this stage the transmitter could
directly transmit $c^{'}$ and the intended receiver would use the
inverse of the RSA system first and then the inverse of the unit
encryption. Alternatively the $c^{'}$ or $c$ could be unit-encrypted using
secret unit  of the transmitter and then transmitted; this would enable
a digital signature to be generated. It is of course also possible to encrypt
$m$ with the secret unit of the transmitter and this would also enable
a digital signature to be generated.

\subsection{Combining group ring cryptographic system and coding.}
It is also possible to combine one of these cryptographic systems and an
(error-correcting) coding system codes in
the one system.

Thus: Given $w$ to be encrypted. Then:  $encode$, $encrypt$,
 $transmit$,  $decrypt$,
 $decode$;  or alternatively $encrypt$, $encode$, $transmit$, $decode$,
 $decrypt$.

There are a number of ways and methods of doing this and it possible
to  use either the group structure or
the ring structure or indeed both for the coding. 

One way is to use the the extra group
elements which are not used in the expression for the data to be encrypted. 
Suppose the word to be encrypted and coded is $w = \di\sum_{i=0}^t\alpha_ig_i$ 
in the group ring of the group $G$ which has $n$ elements where $n > t$. We 
then consider $\hat{w} = w + \di\sum_{j=t+1}^n\beta_jg_j = w + w_1$ (say). 
Then the $w_1$ can be used for coding and error-correcting codes and 
the $\beta_j$ are chosen suitably for the coding. 

Alternatively we can use the ring structure for encoding. Suppose for
example the words we wish to encode (and encrypt) are words in $\Z_t$.
We choose our ring to be $\Z_s$ for some $ s > t$ or the integers
$\Z$ itself. Encode then the
coefficients into $\Z_s$ or into $\Z$.  The process is then: encode, encrypt,
transmit, decrypt and decode.  

Polynomial codes can be combined very nicely with cyclic group rings.
Suppose $G$ is cyclic of order $n$ generated by $g$. Then $RG \cong
R[x]/<x^n -1>$.      
Many of the well known error-correcting codes have been characterised
as ideals in the group ring of certain finite groups.




\paragraph{\bf Order.}
The units can be chosen to have infinite order and 
taking powers of these will not help with finding keys.  If the units
are taken over $\Z_n$ then the 
order will be finite but if $n$ is large this order is also large and
not obtainable by formula.   (In RSA and
other known systems the keys or powers used have finite order although
this order is very large.)


\paragraph{Disguising the size.}
Suppose $u$ is a unit which is to be used for encryption and has
the form $u = \di\sum_{i=1}^{n}\alpha_i g_i$.

It is also possible to disguise $n$ and this can be done as follows:

Choose {\em random numbers} $\beta_j$ for $j=n, \ldots , s$. 
Replace $\alpha_i$ for $i \leq n$ by 
$\alpha_i - \beta_{n+i} - \beta_{n+2i} - \ldots - \beta_{n+ti}$ where
$n+ti \leq s$ and 
$n+(t+1)i > s$. Then use $u = \di\sum_{i=1}^{n}\alpha_ig_i +
\di\sum_{j=n+1}^s\beta_jg_j$, where the elements $g_j$ for $j > n$ are
chosen such that $g_i = g_{n+i} = \ldots = g_{n+it}$.\footnote{This
  could be done for example when $g_i = g^i$ and all that is necessary
to know is that $g^i*g^j = g^{i+j}$ without revealing the order of $g$.
Then the person 
wishing to send the message will simply apply a polynomial
multiplication to their message 
without reducing modulo $g^n = 1$. The message received is reduced
modulo $g^n = 1$ and decrypted.}
The person encrypting the message need not know these relations on
 $g_k$ but would know the multiplication of the $g_k$.
 Then  $n$ can be kept secret. Variations on this theme of adding ``pseudo''
coefficients are also possible.

We are thus able to work in a system without revealing its
order and all we need to know is how the multiplication is performed in
the group ring.  Now \cite{buch} section ``Generalization" on
page 153 discusses this type of difficulty.  

\paragraph{\bf Units and types.}
There are many known units in group rings  and there is extensive literature on
these. Non-commutative systems and units as well as commutative ones 
are available. However given a group ring it is not possible to say
what the units are. Known ones are constructed by formula but
combinations of these are not available by formula. 

Cyclic group rings can be used but also combinations of cyclic
group rings within other group rings.  If $G$ is cyclic of order
$n$ then $RG \cong R[x]/<x^n - 1>$.  Hence if $R$ is a Euclidean Domain
an attack may be made on the cyclic group ring system using the
Extended  Euclidean Algorithm. Thus if a cyclic group ring is used
on its own, the order of the ring and/or the order of the group,
preferably both, should be large just like in RSA the size of the
primes involved must be large.

Units can be
constructed by formula. 
 These can be combined to give others which
cannot be obtained by formula. 
The formula for the inverse could contain additional difficulty, for
example  ``RSA
difficulty'' in the sense that given $n$, the order of the group $G$, 
information on the prime factorisation of $n$, as
well as information on the constructing formula, will be necessary,  
in order to construct the inverse of the unit. Once the inverse of a
set of units
is known then the inverse of the product of these units is known but
a formula may not be used from the product to produce its
inverse. 

In practical cases once the keys have been
calculated the encryption and decryption can be done in $O(n\log n)$ time
at worst; methods exist which improve on this time. In certain cases
existing hardware and software for polynomial multiplication over
$\Z$,  $\Z_n$ or finite fields can be used for fast encryption/decryption.  

\subsection{\bf Variations.}

Many variations and examples are possible which remain within the
concept and scope of the group ring method.

\paragraph{\bf Advantages.}
\begin{itemize}
\item A new method for public key encryption. When the keys have been
  constructed the time for encryption and decryption is short and
  faster than existing methods.
\item The units or cryptographic systems  can be combined to form new
  units  or cryptographic
  systems which are of a different {\it type} to the constituent units
  or cryptographic systems. This has advantages over existing systems.
\item Powers of the units can be used just as easily as the units
  themselves but these powers are more difficult to attack as the
  difficulty of the discrete logarithm problem comes into play.
\item These units or cryptographic systems can be combined with
  existing systems to form new systems which are even more powerful
  than the constituents. Thus combining with RSA will have RSA
  security and the group ring security combined and the system is not
  of either type.
\item There are many units and group rings available of many different
  types which
  can be used. These can be non-commutative as well as
  commutative. Non-commutative systems for public key cryptography do
  not exist at present. Non-commutative systems are more difficult to break.
\item In certain situations the length of system may be hidden. Systems
  of different lengths may be combined.
\end{itemize}

\paragraph{Summary}
\begin{enumerate}
\item A new  public key cryptographic system using group rings. A
  group ring  $RG$ is chosen. A unit of this group ring is generated
  which is then used as the public key. The method of generation of
  the unit or public key is such that the inverse of the unit which is
  then the private key, cannot be obtained from the public key. The
  plaintext to be transmitted is converted into a (unique) group ring element
   and the chosen unit (public key) is used to generate the ciphertext. The
  ciphertext is transmitted and only the person holding the private
  key, which is the inverse of the unit, can decrypt the ciphertext. 
\item Two or more chosen units of the same or of different types may be
  combined to form a new unit which is not of any of the types of the
  constituent units. This new unit is then used as the public key and
  is not of any known type. The private key is kept secret and is the
  combination in a certain order of the inverses of the constituent
  units.
\item The group rings used may be commutative or
  non-commutative. Known units in the non-commutative case may be 
combined together  (or
  with commutative units) to give a new unit not of the same type
  as any  of the constituents. This new unit is
  then the public key and the private key is obtained from the
  inverses of the constituents; the public key would not reveal either the
  inverses of the constituents nor the types of the constituents.

\item Given a unit of a particular type, a power of this unit and
      products of such powers may be
  used as the public key. Using a power of a unit also introduces the
  {\em difficulty of the discrete logarithm problem} so that knowing
  the power of the unit is not enough to know the unit itself. This
  power of  a unit may also be combined with other
  keys or units which have been generated according to this or
  previous methods to give a new unit and hence a new public key and a
  new private key; this new system  has both the difficulty of
  trying to find the inverse of a unit and the difficulty of the discrete
  logarithm problem.

\item A method exists to disguise the size of the units or keys. A
  method exists to deal with large numbers.

\item These units or keys generated may
  be  combined with existing public key cryptosystems such as
  RSA. The new combined systems  have more security than each
  of the constituent systems and will not be of either type.
  Thus, for example, a system combining
  group ring public key and RSA public key will have more security
  than the group ring public key system or the RSA system and will not
  be of either type.
. 
\item In many cases these group ring public key cryptosystems can be
  combined with coding and
  error-correcting codes to give one system with both public key and
  coding all in one system.
\item Many variations and examples are possible.

\end{enumerate} 
Examples are now given  which illustrate the
techniques involved. These should be taken as illustrative and not
exhaustive. Also the lengths and numbers are kept small so they
can be displayed. 
\section{Examples} This examples were constructed with the help of
various Computer Algebra Systems such as GAP or MATLAB. 
The numbers in the examples and
the lengths are kept
reasonable here so that they can be displayed. 

Lines beginning with \begin{verbatim}# \end{verbatim} are comments. Other 
lines are program outputs or inputs.
\begin{itemize}
\item[\bf Example 1]: This is an example of encrypting and decrypting using a 
unit in the group ring.
\begin{verbatim}

#this is the logfile of the running of groupring3, 
#which constructs the objects and does the encryption and decryption
#the time in millesecs to perform calculations is also given
# h is the  public key (encrypting function),
# hinv is private key which is used to decrypt. These have been constructed    
#from a formula in the group ring and involves parameters which are not 
# revealed and do  not need to be revealed.
 
h;
-408-402*g-374*g^2-298*g^3-144*g^4+94*g^5+374*g^6+606*g^7+697*g^8+606*g^9+
374*g^10+94*g^11-144*g^12-298*g^13-374*g^14-402*g^15

hinv;
-13464+5106*g+9622*g^2-12470*g^3-144*g^4+12674*g^5-9622*g^6-5310*g^7+13753*g^8
-5310*g^9-9622*g^10+12674*g^11-144*g^12-12470*g^13+9622*g^14+5106*g^15

#just to check privately that these are inverse of one another
h*hinv;
1
# we encrypt sequence (next)as an example


sequence:= [-3,12,-16,72,-123,1,-1,0,1234,-17,143,0,64,-173,13,-234]
#this could be in any base. It is left in integers.
#write this as a group ring element:

r:= -3 + 12*g -16*g^2 +72*g^3 - 123*g^4 + g^5 - g^6 + 1234*g^8 - 17*g^9 + 
143*g^10 + 64*g^12 - 173*g^13 + 13*g^14 - 234*g^15
-3+12*g-16*g^2+72*g^3-123*g^4+g^5-g^6+1234*g^8-17*g^9+143*g^10+64*g^12-173*g^
13+13*g^14-234*g^15

#encrypt

x:=r*h;

1033182+949413*g+646149*g^2+228128*g^3-179124*g^4-488007*g^5-663825*g^6
-718750*g^7-688787*g^8-614702*g^9-516410*g^10-379628*g^11-164099*g^12+
153119*g^13+534225*g^14+870088*g^15

#the above [1033182,949413,646149, ...,870088] is what would be sent 
#instead of r.

#to decrypt x multiply by hinv
y:= x*hinv ;
-3+12*g-16*g^2+72*g^3-123*g^4+g^5-g^6+1234*g^8-17*g^9+143*g^10+64*g^12-173*g^
13+13*g^14-234*g^15

#check if decryption is done correctly

y=r;
true

#check total runtime; 
Runtime;
10

#answer is in millesecs

#check runtime for encryption and decryption (of sequence)

crypttime;
0

#this took no time at all!

\end{verbatim}
\item[\bf Example 2]
The following example shows how to combine a group ring unit and  RSA method.

\begin{verbatim}
Read("combo");
#Read in the program to perform all calculations
n;
78013681
p;
7459
q;
10459
#n is product of p and q
#Phi(n) = 77995764

e; 
5
#this is exponent which is 5; the inverse of e mod Phi(n) is d = 15599153

#RSA public key is (n,e) = (78013681,5)
#RSA private key is (d) = (15599153)

#Now for the Units which were constructed:
h;
2983+1407*g-573*g^2-2308*g^3-3301*g^4-3301*g^5-2308*g^6-573*g^7+1407*g^8+
2983*g^9+3585*g^10

hinv;
-14659+22389*g-14659*g^2-3190*g^3+18832*g^4-21472*g^5+9295*g^6+9295*g^7
-21472*g^8+18832*g^9-3190*g^10

#Unit public key is h
#Unit private key is hinv
#Combination public key is (n,e,h)
#Combination private key is (d,hinv)
#P = 1231 was encrypted (and decrypted) as an example
#Under RSA it goes to C:= P^e mod n

15134643

#this now goes into group ring format:
r:=
1+5*g+g^2+3*g^3+4*g^4+6*g^5+4*g^6+3*g^7

#This is now encrpypted to (in group ring form) r*h

encrypted;
-11135+11911*g+31358*g^2+41330*g^3+38402*g^4+22982*g^5-201*g^6-23410*g^7
-38792*g^8-41440*g^9-30978*g^10

#or in sequence form:
trans;
[ -11135, 11911, 31358, 41330, 38402, 22982, -201, -23410, -38792, -41440, 
  -30978 ]

#this is transmitted for r = 1231

#to decrypt: First decrypt in group ring:
defirst;
1+5*g+g^2+3*g^3+4*g^4+6*g^5+4*g^6+3*g^7
#In coefficients:
coeff;
[ 1, 5, 1, 3, 4, 6, 4, 3 ]
#and then the number
 w;
15134643
#now decrypt by inverting RSA
#w^d mod n is:
backtostart;
1231

backtostart = P;
true
#this is what we encrypted 

#note time:
Runtime;
20

# The whole thing, including constructions the keys, is done in about 
# 20 millesecs without any special techniques which could have shortened 
# the calculations

\end{verbatim}
\item[\bf Example 3] The next example combines units of different types. 
Then new unit obtained is not of any of these types nor of any other
 known type.
\begin{verbatim}

#We combine a B-unit with a H-unit to get a new unit
#which is used for encryption/decryption. This new unit is
#not of either form i.e is not a B-unit nor a H-unit.
#The order of the group is 16; the order of the group ring is of course 
#infinite. h is the B-unit and hinv is its inverse; hh
#is the H-unit and hhinv is its inverse.
#the public key is then (h*hh) and the private key is
#hhinv*hh. We could take a power of h or hh to make it more difficult and
# we are then within the difficulty of discrete logarithm problem on 
#powers of group ring elements 

Read("groupring3");

h;
-408-402*g-374*g^2-298*g^3-144*g^4+94*g^5+374*g^6+606*g^7+697*g^8+606*g^9+
374*g^10+94*g^11-144*g^12-298*g^13-374*g^14-402*g^15

hinv;
-13464+5106*g+9622*g^2-12470*g^3-144*g^4+12674*g^5-9622*g^6-5310*g^7+13753*g^8
-5310*g^9-9622*g^10+12674*g^11-144*g^12-12470*g^13+9622*g^14+5106*g^15

Read("hunit1");

#the h-unit is hh and its inverse is hhinv

hh;
1-g+g^3-g^4+g^5-g^7+g^8

hhinv;
1-g^2-2*g^3-2*g^4-2*g^5-g^6+g^8+g^9+g^10+g^11+g^12+g^13+g^14+g^15

#combine h and hh to form new unit which is not of either form

newunit:=h*hh;
25+18*g-18*g^2-66*g^3-88*g^4-66*g^5-18*g^6+18*g^7+25*g^8+17*g^9+18*g^10+31*g^
11+39*g^12+31*g^13+18*g^14+17*g^15

#the inverse of this unit is hhinv*hinv; these do not have to
#be multiplied out but could be used one-by-one for decryption
#the private key is then:
newunitinv:=hhinv*hinv;
25-3497*g+2663*g^2+1459*g^3-3791*g^4+1459*g^5+2663*g^6-3497*g^7+25*g^8+3462*g^
9-2663*g^10-1424*g^11+3742*g^12-1424*g^13-2663*g^14+3462*g^15

#check privately:

newunit*newunitinv;
1

#to encrypt the sequence 
#[-12,-234,345,-435,0,165,-142,43,-17,-12,456, -2341,-321,23,-76];
#This sequence could be in any base 
sequence:= [-12,-234,345,-435,0,165,-142,43,-17,-12,456, -2341,-321,23,-76]
#write as group ring element:

-12-12*g-234*g^2+345*g^3-435*g^4+165*g^6-142*g^7+43*g^8-17*g^9-12*g^10+456*g^
11-2341*g^12-321*g^13+23*g^14-76*g^15

#now encrypt
x:= r*newunit;
185871+165276*g+68927*g^2-21364*g^3-51052*g^4-34033*g^5-26102*g^6-48807*g^7
-73742*g^8-67942*g^9-44554*g^10-41339*g^11-63822*g^12-63651*g^13+1671*g^14+
112093*g^15

#thus sequence 
[185871,165276,68927,-21364,-51052,-34033,-26102,-48807,-73742,-67942,-44554,
-41339,-63822,-63651,1671,112093] 
#is transmitted 

#to decrypt multiply by newunitinverse

y:=x*newunitinv;

-12-12*g-234*g^2+345*g^3-435*g^4+165*g^6-142*g^7+43*g^8-17*g^9-12*g^10+456*g^
11-2341*g^12-321*g^13+23*g^14-76*g^15

 #check we got back;

y=r;
true

Runtime;
30
#this is in millesecs
\end{verbatim}
\item[\bf Example 4] These examples shows how to combine group ring 
encryption 
and coding. The process is:  encrypt, encode, transmit (with errors in
transmission),
decode (correcting errors in transmission) and decrypt. The first uses
the ring structure for the coding and the second one uses the group
structure for the coding. Many codes could be used but we have
used fairly simple ones for illustrative purposes.

\medskip 
First example of this type:
\begin{verbatim}

#Will use a Hamming code on each coefficient.
Hamming:= HammingCode(3, GF(2));
a linear [7,4,3]1 Hamming (3,2) code over GF(2)

#The following is to be encrypted, encoded, transmitted, decoded
# and decrypted. We assume errors occur in transmission for illustration.

r:= 11+15*g+12*g^2+8*g^3+13*g^5+11*g^6+7*g^7+13*g^8+4*g^9+7*g^10

# want to encrypt and encode this r. First encrypt using the unit h.

rencrypt:= r*h; #we can take everything mod 16 (at the moment)

11+14*g+5*g^2+7*g^3

# encode all the coefficients of rencrypt using Hamming code
# rencrypt encoded is the following

rencode:= 51 + 22*g + 37*g^2 + 15*g^3;

# so this is transmitted; Some errors may be made on transmitting the binary
# codes  of the numbers: Suppose the following is received

rencodeerror:= 19 + 44*g + 91*g^2 + 79*g^3;

19+44*g+91*g^2+79*g^3

#decode each coeficient
 Decode(Hamming, 19,44,91,79);
11, 14, 5, 7.

# now have the correct version.

corrected:= 11 + 14*g + 5*g^2 + 7*g^3; 
11+14*g+5*g^2+7*g^3

# now decrypt
#first note all coefficients are now less than 16

decrypt(corrected);

11+15*g+12*g^2+8*g^3+13*g^5+11*g^6+7*g^7+13*g^8+4*g^9+7*g^10

#we are back
Runtime;
10
#10 millesecs to do it all.
\end{verbatim}
Second example of this type:
\begin{verbatim}
#encrypt [1,0,1,1], encode, transmit, decode, decrypt
start:= 1 + g^2 + g^3;
#encode with unit h
r:= start*h ;
Z(2)^0+g+g^2+g^4+g^6+g^7+g^10
# now encode with Hamming [15,13] code over GF(2)
 C2:= HammingCode(4, GF(2));
a linear [15,11,3]1 Hamming (4,2) code over GF(2)
#as a sequence r is
seq:= [1,1,1,0,1,0,1,1,0,0,1];
#encode this with the Hamming code

seqencode:= seq*C2;
          
[ 1 0 1 0 1 1 0 0 1 0 1 1 0 0 1 ]

#this is transmitted but an error  occurs in the first position
#The following then is what is received.

seqcodeerror:=          
[0, 0, 1, 0, 1, 1, 0, 0, 1, 0, 1, 1, 0, 0, 1 ]

#now decode this:
Decode(C2, seqcoderror);

          
[ 1 0 1 0 1 1 0 0 1 0 1 1 0 0 1 ]

#error has been corrected
# now decrypt:
#  
decrypt:= 1 + g^2 + g^3

#check if we got back to where we started:

decrypt = start;
true

\end{verbatim}
\item[\bf Example 5] This example takes a large power of a unit and
  uses this an the public key. This introduces the difficulty of the
  discrete logarithm problem.
\begin{verbatim}
#This is an example where a unit is constructed and then a power if this unit
#is constructed. The numbers and expressions are too large printed out here.
# n is the size of the group
# d is the length of the original unit hh; hhinverse is the inverse of 
# hh. Then hh^127 is taken; q = 127 is the power of hh worked out
# hhpower = hh^127; this is the public key; too large to print out
# hhpowerinverse is the inverse of hhpower which is of course 
# equal to hhinverse^127.
# then hhpowerinverse is the private key
# as an example we take the sequence [1,1,1...] where 1 appears 4096 times
# and encrypt and decrypt this.
# no sophisticted time saving or hard-coding was used.
n;
4096
d;
511
q;
127
Read("alunit22");
#This reads in the programme to do all the calculations.
# r is constructed as r:=Sum([0..(n-1)],t->g^t);;
#encrypt
y:= r*hhpower;;
x:= y*hhinversepower;
#check if encryption/decryption worked out
x=r;
true
#as already mentioned the numbers and sequences 
#are not printed as they are too big.
Runtime;
2500
# this is in millesecs; the whole thing took 2.5 sec
\end{verbatim}
\item[\bf Example 6] This example computes non-commutative units and
  combines these to get a public key system.
\begin{verbatim}
 n;
10
#Symmetric group on 10 letters and of order Factorial(n) = 3628800
# is constructed
#Its group ring is then constructed; two noncommutative
#units uab and bau are then constructed  together with their inverses
# uabinv, bauinv 
# Also a cyclic unit h and its inverse hinv are constructed.
# Read in the programme to do the calculations
Read("noncommunit");

uab;;
bau;;

#these uab and bau are too long to be displayed here
# check if they commute
uab*bau = bau* uab;

false

# form a new unit by multiplying these together; its inverse is
#uabinv*bauinv which is *not* the same as bauinv*uabinv
bau*uab;;
# it's too long to display here as is:

uab*bau;;

#we now form new encrypting and decrypting functions by various 
#multiplications of uab, bau, h and their inverse: 

encrypter1:=(uab*bau);;
decrypter:=(bauinv*uabinv);;
encrypter2:=(bau*uab);;
decrypter2:=(uabinv*bauinv);;

#check if these are the same
encrypter1=encrypter2;
false

#check, privately, that these are inverse to one another.
encrypter1*decrypter1;

identity

#identity is the identity element of the group ring

encrypter2*decrypter2;

identity

#take powers: 
encrypter3:=encrypter1^5;;

decrypter3:=decrypter1^5;;

enceypter3*decrypter3;
identity 
# as expected
#multiply all three to get another system:

encrypter4:=uab*bau*h;;
decrypter4:=hinv*bauinv*uabinv;;

# these are all pretty large so are not printed out

encrypter4*decrypter4;
identity 

#as expected

#to show that order is important
defalse:= bauinv*uabinv*hinv;;
# see if this is the inverse of encrypter4
 defalse=decoder4;
false

#How long did all this take? No time-saving techniques or hard-coding
#was used
Runtime;
35
#it took about 35 millesecs. The lengths were in the thousands with the powers.
\end{verbatim}

\end{itemize}

\end{document}